\documentstyle[11pt,leqno]{article}
\newtheorem{theorem}{{\sc Theorem}}
\newcommand{\bt}{\begin{theorem}}
\newcommand{\et}{\end{theorem}}
\newcommand{\newsection}[1]{\setcounter{equation}{0} \setcounter{theorem}{0}
\section{#1}}

\newcommand{\NI}{\noindent}
\newcommand{\bea}{\begin{eqnarray}}
\newcommand{\eea}{\end{eqnarray}}

\def \spec#1 {\mathop{#1}}

\def \b #1 {\bf #1}
\newcommand {\nnb }{\nonumber}

\newcommand {\CC}{\centerline}

\newcommand{\clf}{{\cal F}}

\newcommand{\raro}{\rightarrow}

\newcommand{\vsp}{\vskip 1em}
\newcommand{\vspp}{\vskip 2em}

\newcommand{\be}{\begin{equation}}
\newcommand{\ee}{\end{equation}}
\newcommand{\ben}{\begin{eqnarray*}}
\newcommand{\een}{\end{eqnarray*}}

\begin{document}

\sloppy
\CC {\Large{\bf Nonparametric Estimation of Trend for SDEs with Delay}}
 \CC {\Large {\bf Driven by Fractional Brownian Motion with Small Noise}}
\vsp
\CC {\bf B.L.S. PRAKASA RAO}
\CC {\bf CR Rao Advanced Institute of Research in Mathematics,}
\CC{\bf  Statistics and Computer Science, Hyderabad 500046, India}
\vsp
\NI{\bf Abstract:} We investigate the problem of nonparametric estimation of the trend for stochastic differential equations with delay and driven by a fractonal Brownian motion through the method of kernel-type estimation for the estimation of a probability density function.
\vsp
\NI{\bf Keywords and phrases}: Nonparametric estimation; Estimation of trend; Stochastic differential equation with delay; Kernel method of estimation; Fractional Brownian motion. 
\vspp
\NI {\bf AMS Subject classification (2010): } Primary 62M09, Secondary 60G22.
\vsp
\newsection{Introduction}
Guschin and Kuchler (1999) investigated asymptotic inference for linear stochastic differential equations (SDEs) with a time delay of the type
$$dX(t)= (aX(t)+bX(t-1))dt +dW_t, t \geq 0$$
driven by the standard Brownian motion $W= \{W_t, t\geq 0\}$ with the initial condition $X(t)=X_0(t), -1\leq t \leq 0$ where  $\{X_0(t), t \in R\}$ is a continuous process independent of the process $W.$ They investigated  the asymptotic properties of the maximum likelihood estimator (MLE) of the parameter $\theta=(a,b)$ based on the observation of the process $\{X(t), 0\leq t \leq T\}.$ They showed that the asymptotic behaviour of the MLE depends on the ranges of the values of $a$ and $b.$ Prakasa Rao (2008) considered the problem of estimation of the parameter (a,b) for linear SDE of the form
\be
dX(t)= (aX(t)+bX(t-1))dt +dW^H_t, t \geq 0
\ee
driven by the fractional Brownian motion $W^H= \{W^H_t, t\geq 0\}$ with known Hurst index $H$ and with the initial condition $X(t)=X_0(t), -1\leq t \leq 0$ where  $\{X_0(t), t \in R\}$ is a continuous process independent of the process $W^H.$ Asymptotic properties of the MLE for the parameter $\theta = (a,b)$ is studied in Prakasa Rao (2008). For some details, see Prakasa Rao (2010), pp.156-163. Applying the results in Mohammed and Scheutzow (1990), it can be shown that there exists a unique solution $X=\{X(t), t \geq -1\}$ of the equation (1.1) and it can be represented in the form
\be
X(t)=x_0(t)X_0(0)+b\int_{-1}^0x_0(t-s-1)X_0(s)ds+\int_0^tx_0(t-s)dW_s^H, t\geq 0.
\ee  
This process has continuous sample paths for $t\geq 0$ almost surely and conditionally on $X_0,$ the process $X$ is a Gaussian process. Furthermore the function $x_0(t)$ defined for $t \geq -1$ is the fundamental solution of the differential equation 
\be
\frac{dx(t)}{dt}=ax(t)+bx(t-1), t >0 
\ee
subject to the condition $x(0)=1, x(t)=0, t \in [-1,0].$ Guschin and Kuchler (1999) gives a complete discussion on the existence and representation of the fundamental solution of the equation (1.3). 
\vsp
Kutoyants (2019) studied nonparametric estimation for SDE with delay driven by a standard Wiener process with small noise. We wlll now investigate a similar problem  for SDE with delay driven by a fractional Brownian motion  with small noise.  
\vsp
\newsection{Preliminaries} Consider the SDE with delay given by
\be
dX_t=S(t,X_{t-\tau})dt+\epsilon \; dW^H_t, X_s=x_0, s \leq 0, 0\leq t \leq T 
\ee
where $W^H$ is a fractional Brownian motion with known Hurst index $H$ and the function $S(t,x)$ is unknown. We assume that the delay $\tau \geq 0$ is known. Consider the ordinary differential equation 
\be
\frac{dx_t}{dt}= S(t,x_{t-\tau}), x_s=x_0, s \leq 0, 0\leq t \leq T. 
\ee
For any given value $x,$ let $t=t_x$ be defined by the equation $x_{t_x-\tau}=x.$ The problem is to estimate the function $S(t_x,x)$ based on observations $X^T=\{X_t,0\leq t \leq T\}.$ Observe that the equation (2.1) can be represented in the form
\be
X_{t-\tau}= x_0+\int_0^{t-\tau}S(s,X_{s-\tau})ds+\epsilon W^H_{t-\tau}, t \geq \tau.
\ee
Parametric inference for SDEs of delay type and driven by Wiener process was investigated by several authors. A survey of these results is given in Kutoyants (2005). Nonparametric inference for such processes with small noise was surveyed in Kutoyants (1994, 2004). Parametric inference for SDEs of delay type and driven by fractional Brownian motion was investigated by Prakasa Rao (2008) and a discussion was presented in Prakasa Rao (2010). We now initiate the study of nonparametric inference for processes with time delay but driven by a fractional Brownian motion.
\vsp
\newsection{Main Results}
We assume that the  the following conditions hold.\\
\NI $(C_1)$ The function $S(t,x)$ is positive and $\alpha \equiv \inf_{t \in [0,T],x}S(t,x) >0.$\\
\NI $(C_2)$ The function $S(t,x)$ has two continuous bounded derivatives with respect to $t$ and $x.$ \\
\NI $(C_3)$ The delay parameter $\tau$ and the Hurst index $H$ are known.\\
\vsp
Note that the condition $(C_2)$ implies that there exists constants $M>0$ and $L>0$ such that 
\be
|S_t^\prime(t,x)| \leq M; \;\;\mbox{and}\;\;|S_x^\prime(t,x)|\leq L
\ee
where $S_t^\prime(t,x)$ denotes the partial derivative of $S(t,x)$ with respect to $t$ and $S_x^\prime(t,x)$ denotes the partial derivative of $S(t,x)$ with respect to $x.$ We now consider the problem of estimation of the function $S(t,x)$ for $t >\tau.$ For $t \leq \tau,$ the function $S(t,x_t)=S(t,x_0)$ and hence a function of one variable $t \in [0,\tau].$ This problem was investigated in Mishra and Prakasa Rao (2011) and in Section 5.3 in Prakasa Rao (2010). Here after, we assume that $t> \tau.$
\vsp
Let $G(.)$ be a real-valued function which is bounded and satisfies the conditions:\\

\NI(i) $\int_{-1}^1 G(u) du=1,$\\
\NI(ii) $\int_{-1}^1 u \;G(u) du=0,\;\;\mbox{and}$\\
\NI(iii) $G(u)=0, |u|\geq 1.$

Let $\phi_\epsilon >0 $ be a bandwidth function tending to zero as $\epsilon \raro 0.$ We now construct a kernel type estimator for the function $S(t,x)$ for $t >\tau$ based on the kernel $G(.)$. For a study of kernel-type estimators and their properties for the estimation of a probability density function, see Prakasa Rao (1983). Let 
\be
t_{x,\epsilon}= \{t>\tau: X_{t-\tau}\geq x\}.
\ee
If $\sup_{0\leq s\leq T-\tau} X_s < x,$ then, define $t_{x,\epsilon}=T-\tau.$ Define the kernel-type estimator
\be
S_\epsilon(t,x)=\frac{1}{\phi_\epsilon}\int_0^TG(\frac{s-t_{x,\epsilon}}{\phi_\epsilon})dX_s
\ee
as an estimator for the function $S(t,x)$ for $t>\tau.$ We will prove a lemma which will be used in the sequel. In the following discussion, the constant $C$ may be different from one inequality to another.
\vsp
\NI{\bf Lemma 3.1:} Suppose that the function $S(t,x), 0\leq t \leq T, x \in R$ has a continuous bounded derivative with respect to $x$ and the the process $\{X_t, 0\leq t \leq T\}$ is a solution of the SDE given by  (2.1). Let $x_t,0\leq t \leq T$ be a solution of the differential equation (2.2). Then there exists a constant $C>0,$ such that 
\begin{description}
\item (i) $|X_t-x_t| \leq C \epsilon \sup_{0\leq s \leq T}|W^H_s|, 0\leq t \leq T$\\
\item (ii) $E|X_t-x_t|^2 , \leq C \epsilon^2, 0\leq t \leq T$\\
\end{description}
\vsp
\NI{\bf Proof:} It is easy to see that
\bea
|X_t-x_t| &\leq & \int_0^t|S(s, X_{s-\tau})-S(s,x_{s-\tau})|ds+ \epsilon |W^H_t|\\\nnb
&\leq & L \int_0^t |X_{s-\tau}-x_{s-\tau}|ds + \epsilon \sup_{0\leq s \leq T}|W^H_s|\\\nnb
&\leq & L \int_0^t |X_s-x_s|ds + \epsilon \sup_{0\leq s \leq T}|W^H_s|\\\nnb 
\eea
almost surely. Applying the Grownwall lemma, it follows that there exists a constant $C>0$ such that 
$$|X_t-x_t| \leq C \epsilon \sup_{0\leq s \leq T}|W^H_s|, 0\leq t \leq T$$
which in turn implies that 
\bea
E|X_t-x_t|^2 & \leq & C \epsilon^2 E[\sup_{0\leq s \leq T}|W^H_s|^2], 0\leq t \leq T\\\nnb
& \leq & C \epsilon^2 T^{2H}, 0\leq t \leq T\\\nnb
\eea
by the estimate for the supremum of a fractional Brownian motion over an interval $[0,T]$(cf. Novikov and Valkeila (1999); Prakasa Rao (2010), Prop. 1.9). 
\vsp
\NI{\bf Theorem 3.1:} Suppose the conditions $(C_1)-(C_3)$ hold. Let $\phi_\epsilon= \epsilon ^{1/(3-H)}$ and $\hat x= x_0+(T-\tau)\alpha.$ Then,for any $x \in (x_0, \hat x),$ the estimator $S_\epsilon(t,x)$ is $L_2$-consistent and
\be
E[|S_\epsilon(t,x)-S(t,x)|^2]\leq C \epsilon^{4/(3-H)}
\ee 
for some positive constant $C.$
\vsp
\NI{\bf Proof:} Define the set $\Gamma_\epsilon= \{ \omega: \phi_\epsilon^{-1}t_{x,\epsilon}>2, \phi_\epsilon^{-1}(T-t_{x,\epsilon})>2 \}.$ For $\omega \in \Gamma_\epsilon$, we now estimate the integral 
$$\frac{1}{\phi_\epsilon}\int_0^TG(\frac{s-t_{x,\epsilon}}{\phi_\epsilon})S(s,x_{s-\tau})ds.$$
Applying the transformation $s=t_{x,\epsilon}+\phi_\epsilon u$ and using the Taylor expansion,  we get that
\bea
\lefteqn{\frac{1}{\phi_\epsilon}\int_0^TG(\frac{s-t_{x,\epsilon}}{\phi_\epsilon})S(s,x_{s-\tau})ds}\\\nnb
&=& \int_{-\frac{t_{x,\epsilon}}{\phi_\epsilon}}^{\frac{T-t_{x,\epsilon}}{\phi_\epsilon}}  G(u) S(t_{x,\epsilon}+\phi_\epsilon u, x_{t_{x,\epsilon}+\phi_\epsilon u-\tau})du\\\nnb
&=& \int_{-1}^1G(u) S(t_{x,\epsilon}+\phi_\epsilon u, x_{t_{x,\epsilon}+\phi_\epsilon u-\tau})du\\\nnb
&=& S(t_{x,\epsilon},x_{t_{x,\epsilon}-\tau})\\\nnb
&&\;\; +\phi_\epsilon \int_{-1}^1 u\;G(u)du \;\;[S^\prime_t(t_{x,\epsilon},x_{t_{x,\epsilon}-\tau}})+ S^\prime_x(t_{x,\epsilon},x_{t_{x,\epsilon}-\tau})S(t_{x,\epsilon},x_{t_{x,\epsilon}-2\tau)]\\\nnb
&&\;\;\;+ \phi_\epsilon^2 R^{(1)}_\epsilon(t_{x,\epsilon})\\\nnb
&=& S(t_{x,\epsilon},x_{t_{x,\epsilon}-\tau})+ \phi_\epsilon^2 R^{(1)}_\epsilon(t_{x,\epsilon}),\\\nnb
\eea  
where $R^{(1)}_\epsilon(t_{x,\epsilon})$ is the remainder term, by the properties of the kernel function $G(.)$ and by the Taylor expansion of the function $S(t_{x,\epsilon}+\phi_\epsilon u, x_{t_{x,\epsilon}+\phi_\epsilon u-\tau})$ by the powers of $\phi_\epsilon.$ Observe that
\be
x=x_{t_x-\tau}=x_0+\int_0^{t_x-\tau}S(s,x_{s-\tau})ds 
\ee
and
\be
x=X_{t_{x,\epsilon}-\tau}=x_0+\int_0^{t_{x,\epsilon}-\tau}S(s, X_{s-\tau})ds + \epsilon W^H_{t_{x,\epsilon}-\tau} 
\ee
Subtracting the terms on the left side  of the equation (3.8) from the left side of the equation (3.9) and equating to the difference of the terms on the right side of the equations (3.8) and (3.9) and rearranging, we obtain that
\be
\int_{t_x-\tau}^{t_{x,\epsilon}-\tau}S(s,x_{s-\tau})ds=\int_0^{t_{x,\epsilon}-\tau}[S(s,X_{s-\tau})-S(s,x_{s-\tau})]ds+\epsilon W^H_{t_{x,\epsilon}-\tau}
\ee
Applying the condition $(C_1)$, it follows that
\bea
|\int_{t_x-\tau}^{t_{x,\epsilon}}S(s,x_{s-\tau})ds|&\leq &\int_0^{t_{x,\epsilon}-\tau}|S(s,X_{s-\tau})-S(s,x_{s-\tau})|ds+\epsilon |W^H_{t_{x,\epsilon}-\tau}|\\\nnb
& \leq & L\int_0^{T-\tau}|X_{s-\tau}-x_{s-\tau}|ds+\epsilon |W^H_{t_{x,\epsilon}-\tau}| \\\nnb
& \leq & \epsilon H(t_{x,\epsilon})\\\nnb
\eea
where the random variable $H(t_{x,\epsilon})$ has bounded moments by  Lemma 2.1 and the maximal inequalities for the  fBm. Note that
\be
|\int_{t_x-\tau}^{t_{x,\epsilon}-\tau}S(s,x_{s-\tau})ds|\geq \alpha |t_{x,\epsilon}-t_x|
\ee
by the conditions on the function $S(s,x)$. Applying the inequality derived above, it follows that
\be
|t_{x,\epsilon}-t_x|\leq \epsilon \alpha^{-1}H(t_{x,\epsilon}) 
\ee
and the last term tends to zero as $\epsilon \raro 0.$ Observe that
\bea
S_\epsilon(t,x) &= &\frac{1}{\phi_\epsilon}\int_0^T G(\frac{s-t_{x,\epsilon}}{\phi_\epsilon})S(s,x_{s-\tau})ds\\\nnb
&&\;\;\;\;+\frac{\epsilon}{\phi_\epsilon}\int_0^TG(\frac{s-t_{x,\epsilon}}{\phi_\epsilon})dW^H_s\\\nnb
\eea
and hence, on the set $\Gamma_\epsilon,$ it follows that
\bea
\;\;\;\;\\\nnb
|S_\epsilon(t,x)-S(t,x) &=& |S_\epsilon(t_x,x_{t_x-\tau})-S(t_x,x_{t_x-\tau})|\\\nnb
&\leq & |\frac{1}{\phi_\epsilon}\int_0^T G(\frac{s-t_{x,\epsilon}}{\phi_\epsilon})[S(s,x_{s-\tau})-S(t_x,x_{t_x-\tau})]ds|\\\nnb
&&\;\;\;\;+|\frac{1}{\phi_\epsilon}\int_0^T G(\frac{s-t_{x,\epsilon}}{\phi_\epsilon})[S(s,X_{s-\tau})-S(s,x_{s-\tau})]ds|\\\nnb
&&\;\;\;\;\;+ |\frac{\epsilon}{\phi_\epsilon}\int_0^TG(\frac{s-t_{x,\epsilon}}{\phi_\epsilon})dW^H_s|\\\nnb
&\leq & |S(t_{x,\epsilon},x_{t_{x,\epsilon}-\tau})-S(t_x,x_{t_x-\tau})| +\phi_\epsilon^2 |R^{(1)}_\epsilon|+|R^{(2)}_\epsilon|+ |R^{(3)}_\epsilon|\;\;\mbox{(say)}\\\nnb
\eea
from the inequality (3.7). Let $\chi_{\Gamma_\epsilon}$ denote the indicator of the set $\Gamma_\epsilon.$ Observe that
\ben
\lefteqn{E[\chi_{\Gamma_\epsilon}|S(t_{x,\epsilon},x_{t_{x,\epsilon}-\tau})-S(t_x,x_{t_x-\tau})|^2]}\\\nnb
&\leq & 2ME[\chi_{\Gamma_\epsilon}|t_{x,\epsilon}-t_x|^2]+ 2L E[|x_{t_{x,\epsilon}-\tau}-x_{t_x-\tau}|^2]\\\nnb
&\leq & C\epsilon^2 \\\nnb
\een
by the conditions $(C_1)$ and by the inequality given in (3.13). Furthermore, there exists an absolute constant $c(2,H)$ depending on $H$ such that
\bea
\;\;\;\\\nnb
E[\chi_{\Gamma_\epsilon}(R^{(3)}_\epsilon)^2] &\leq & E[(R^{(3)}_\epsilon)^2]\\\nnb
&\leq & \frac{\epsilon^2}{\phi_\epsilon^2} c(2,H) (E[(\int_0^T |G(\frac{s-t_{x,\epsilon}}{\phi_\epsilon})|^{1/H}ds)])^{2H}\\\nnb
& &\;\;\;\;\mbox{(by Memin et al. (2001), Prakasa Rao (2010), Theorem 1.7)})\\\nnb
&=& \frac{\epsilon^2}{\phi_\epsilon^2} \phi_\epsilon^{2H}[\int_{-t_{x,\epsilon}/\phi_\epsilon}^{(T-t_{x,\epsilon})/\phi_\epsilon}|G(u)|^{1/H}du]^{2H}\\\nnb
&\leq & \frac{\epsilon^2}{\phi_\epsilon^{2-2H}}(\int_{-1}^1 |G(u)|^{1/H}du)^{2H}\\\nnb
&\leq & C \frac{\epsilon^2}{\phi_\epsilon^{2-2H}}\\\nnb
\eea
for some constant $C$ depending on $H.$  Let us denote $t_{x,\epsilon}+\phi_\epsilon u$ by $t_{x,\epsilon,u}$ in the following computations. Note that

\bea
E[\chi_{\Gamma_\epsilon}(R^{(2)}_\epsilon)^2] &\leq & E[(R^{(2)}_\epsilon)^2]\\\nnb
&=& |\frac{1}{\phi_\epsilon}\int_0^T G(\frac{s-t_{x,\epsilon}}{\phi_\epsilon})[S(s,x_{s-\tau})-S(s,x_{s-\tau})]ds|^2\\\nnb
&\leq & 2 E\int_{-1}^1 G^2(u) |S(t_{x,\epsilon,u},X_{t_{x,\epsilon, u}-\tau})-S(t_{x,\epsilon,u},x_{t_{x,\epsilon, u}-\tau})|^2du\\\nnb
&\leq & 2 L^2\int_{-1}^1 G^2(u) E|X_{t_{x,\epsilon, u}-\tau}-x_{t_{x,\epsilon, u}-\tau}|^2\\\nnb
&\leq & C\epsilon^2 \\\nnb
\eea
for some positive constant $C$ by an application of Lemma 2.1. Combining the bounds obtained above, we get that
\be
E[\chi_{\Gamma_\epsilon}|S_\epsilon(t,x)-S(t,x)|^2 ]\leq C \phi_\epsilon^4 + C \epsilon^2 + C \frac{\epsilon^2}{\phi_\epsilon^{2-2H}}.
\ee
We now estimate the $P(\Gamma_\epsilon^c).$ Observe that
\bea
P(\Gamma_\epsilon^c) & \leq & P(\phi_\epsilon^{-1} t_{x,\epsilon}\leq 2)+ P(\phi_\epsilon^{-1}(T-t_{x,\epsilon})\leq 2 )\\\nnb
&=& P(t_x-t_{x,\epsilon} \geq t_x-2\phi_\epsilon)+ P(T-t_{x,\epsilon}\leq 2 \phi_\epsilon)\\\nnb
\eea
and hence, for $\phi_\epsilon <\frac{1}{4}\tau,$ it follows that
\bea
P(t_x-t_{x,\epsilon} \geq t_x-2\phi_\epsilon) &\leq & P(|t_x-t_{x,\epsilon}| \geq t_x -\frac{\tau}{2})\\\nnb
&\leq & P(|t_x-t_{x,\epsilon}| \geq \frac{\tau}{2})\\\nnb
&\leq & P(C\epsilon \sup_{0\leq s \leq T}|W^H_s|\geq \frac{\tau}{2})\\\nnb
&\leq & C e^{-c\epsilon^{-2}}\\\nnb
\eea
by Theorem 1.4 in Prakasa Rao (2014) for the  fractional Brownian motion. A similar estimate can be obtained for the term $P(T-t_{x,\epsilon}\leq 2 \phi_\epsilon).$ We can now estimate the mean squared error on the set $\Gamma_\epsilon^c.$ Observe that
\ben
(E[\chi_{\Gamma_\epsilon^c}|S_\epsilon(t,x)-S(t,x)|^2)^2 &\leq & E[\chi_{\Gamma_\epsilon^c}]E|S_\epsilon(t,x)-S(t,x)|^4]\\\nnb
&\leq & C e^{-c \epsilon^{-2}}\\\nnb
\een
The last inequality follows by observing that 
\ben
E|S_\epsilon(t,x)-S(t,x)|^4] \leq 2E|S_\epsilon(t,x)|^4+ 2 [S(t,x)]^4<
\een
and noting that the random variable $S_\epsilon(t,x)$ has polynomial moments of all orders from the conditions on the kernel $G(.)$ and from the bounds on Wiener integrals with respect to the fractional Brownian motion due to Memin et al. (2001). Furthermore the function $S(t.x)$ is bounded by assumption. Combining the inequalities obtained above, we get that
\be
E|S_\epsilon(t,x)-S(t,x)|^2]\leq C\phi_\epsilon^4 + C\epsilon^2 + C \frac{\epsilon^2}{\phi_\epsilon^{2-2H}} + Ce^{-(c/2) \epsilon^{-2}}
\ee
for some positive constants $C$ and $c.$ The optimal rate can be obtained by choosing $\phi_\epsilon $ such that
$$\phi_\epsilon^4= \frac{\epsilon^2}{\phi_\epsilon^{2-2H}}$$
or equivalently $\phi_\epsilon= \epsilon^{1/(3-H)}.$ Hence we obtain the bound
\be
E|S_\epsilon(t,x)-S(t,x)|^2\leq C\epsilon^{4/(3-H)}.
\ee 
\newsection{Comments :}
\NI{\bf Remarks 4.1:}\\
\NI{$(C_4)$} Suppose that the function $S(t,x)$ is Lipschitz and the function $g(t)= S(t,x_{t-\tau})$ has $k$ continuous derivatives and the $k$-th derivative $g^{(k)}(t)$ satisfies 
\be
|g^{(k)}(t)-g^{(k)}(u)|\leq C|u-v|^\beta, u,v \in R 
\ee
for some positive constant $C_0$ and for some $0\leq \beta \leq 1.$ Let $\clf_{k,\beta}$ be the class of functions $g(.)$ satisfying the conditions stated above.\\ 

\NI{$(C_5)$} Suppose that the kernel $G(.)$ satisfies the condition 
$$\int_{-1}^1u^jG(u)du=0, j=1,\dots,k$$ 
in addition to the conditions on $G(.)$ stated earlier.\\ 

Following the arguments stated in Kutoyants (2019), it can be shown that
\be
E[|S_\epsilon(t,x)-S(t,x)|^2]\leq C \phi_\epsilon^{2(k+\beta)}+C \epsilon^2+C\frac{\epsilon^2}{\phi_\epsilon^{2-2H}}
\ee
for some positive constant $C.$ Choosing $\phi_\epsilon$ such that
$$\phi_\epsilon^{2k+\beta}= \frac{\epsilon^2}{\phi_\epsilon^{2-2H}}$$
equivalently
$$\phi_\epsilon= \epsilon^{1/(k-H+\beta+1)},$$
we obtain the inequality
\be
E[|S_\epsilon(t,x)-S(t,x)|^2]\leq C \epsilon^{2(k+\beta)/(k-H+\beta+1)}.
\ee
By a similar analysis, it can be shown that, for any $p>0,$ 
\be
\sup_{S(.,.)\in \clf_(k,\beta)} \epsilon^{-p(k+\beta)/(k-H+\beta+1)}E[|S_\epsilon(t,x)-S(t,x)|^p]\leq C 
\ee
for some positive constant $C.$ It can be shown that , for $t=t_x$ such that $x_{t_x-\tau}=x,$ following the arguments given earlier,
\be
\liminf_{\epsilon \raro 0} \inf_{\hat S_\epsilon(t,x)}\sup_{S(.)\in clf_{k,\beta}}\epsilon^{-2(k+\beta)/(k-H+\beta+1)}E[|\hat S_\epsilon(t,x)-S(t,x)|^2] >0
\ee
implying that the estimator $S_\epsilon(t,x)$ is asymptotically efficient in the sense of Hajek-Lecam (cf. Ibragimov and Khasminskii(1981)).
\vsp
\NI{\bf Remarks 4.2:} Suppose the delay $\tau$ is {\it not known} and we want to estimate the function $f(t)=S(t,x_t-\tau), 0\leq t \leq T.$Suppose the function $f(t)$ is $k$ times continuously differentiable and the $k$-th derivative satisfies the condition (4.1). We can estimate the function $f(t)$ by the kernel estimator
\bea
\hat f_\epsilon(t) &= &\frac{1}{\phi_\epsilon}\int_0^TG(\frac{s-t}{\phi_\epsilon})dX_s\\\nnb
&=& \frac{1}{\phi_\epsilon}\int_0^TG(\frac{s-t}{\phi_\epsilon})S(s,X_{s-\tau})ds +\frac{\epsilon}{\phi_\epsilon}\int_0^TG(\frac{s-t}{\phi_\epsilon})dW^H_s.\\\nnb
\eea
Note that 
$$S(s,X_{s-\tau})=S(s,x_{s-\tau})+O_p(\epsilon)$$
by Lemma 2.1 and 
\be
\frac{\epsilon}{\phi_\epsilon}\int_0^TG(\frac{s-t}{\phi_\epsilon})dW^H_s\leq C \frac{\epsilon^2}{\phi_\epsilon^{2-2H}}
\ee
for some positive constant $C$ from the earlier computations. It is sufficient to study the limit behaviour of the term
\be
J_\epsilon=\frac{1}{\phi_\epsilon}\int_0^TG(\frac{s-t}{\phi_\epsilon})S(s,X_{s-\tau})ds.
\ee
Suppose the kernel $G(.)$ satisfies the condition $(C_5)$ and the conditions stated earlier. Applying the Taylor's expansion, it can be shown that 
\bea
J_\epsilon &=& \int_{-1}^1 G(u)[ f(t)+ \sum_{j=1}^k\frac{\phi_\epsilon^ju^j}{j!}f^{(j)}(t)]du + \phi_\epsilon^kR_\epsilon\\\nnb
&=& f(t)+ \phi_\epsilon^kR_\epsilon\\\nnb
\eea
where $ \phi_\epsilon^k|R_\epsilon|\leq C \phi_\epsilon^{k+\beta}.$ The normalizing function $\phi_\epsilon= \epsilon^{1/(k-H+\beta+1)}$ will give the bound for the mean square error as given in (4.3).
\vsp
\NI{\bf Acknowledgment:} This work was supported under the scheme ``INSA Senior Scientist" of the Indian National Science Academy at the CR Rao Advanced Institute of Mathematics, Statistics and Computer Science, Hyderabad 500046, India.
\vsp
\NI {\bf References}
\begin{description}

\item Guschin, A.A. and Kuchler, U. (1999) Asymptotic inference for a linear stochastic differential equation with time selay, {\it Bernoulli}, {\bf 5}, 1059-1098.

\item Ibragimov, I.A. and Khasminskii, R.Z. (1981) {\it Statistical Estimation: Asymptotic Theory}, Springer, New York.

\item Kutoyants, Yu.A. (1994) {\it Identification of Dynamical Systems with Small Noise}, Kluwer, Dordrecht.

\item Kutoyants, Yu.A. (2004) {\it Statistical Inference for Ergodic Diffusion Processes}, Springer, London.

\item Kutpyants, Yu.A. (2005) On delay estimation for stochastic differential equations, {\it Stochastics and Dynamics},  {\bf 5}, 333-342. 

\item Kutoyants, Yu.A. (2019) On nonparametric estimation for SDE with delay, {\it Publ. Inst. Statist. Univ. Paris}, {\bf 63}, No.2-3, 11-20.

\item Memin, J., Mishura, Y., and Valkeila, E. (2001) Inequalities for the moments of Wiener integrals with respect to a fractional Brownian motion, {\it Statist. Probab. Lett.}, {\bf 51}, 197-206.

\item Mishra, M.N. and Prakasa Rao, B.L.S. (2011) Nonparametric estimation of trend for stochastic differential equations driven by fractional Brownian motion, {\it Statist. Infer. Stoch. Proc.}, {\bf 14}, 101-109.

\item Mohammed, S.E. and Scheutzow, M.K.R. (1990) Lyapunov exponents and stationary solutions for affine stochastic delay equations, {\it Stochastics Stochastic Rep.}, {\bf 29}, 259-283.

\item Novikov, A.A. and Valkeila, E. (1999) On some maximal inequalities for fractional Brownian motion, {\it Statist. Probab. Lett.}, {\bf 44}, 47-54.

\item Prakasa Rao, B.L.S. (1983) {\it Nonparametric Functional Estimation}, Academic Press, Orlando.

\item Prakasa Rao, B.L.S. (2008) Parameter estimation for linear stochastic delay differential equations driven by fractional Brownian motion, {\it Random Oper. Stoch. Equations}, {\bf 16}, 27-38.

\item Prakasa Rao, B.L.S. (2010) {\it Statistical Inference for Fractional Diffusion Processes}, Wiley, London.

\item Prakasa Rao, B.L.S. (2014) Maximal inequalities for fractional Brownian motion, {\it Stoch. Anal. Appl.}, {\bf 31}, 785-799.

\end{description}
\end{document}